\documentclass[11pt]{article}
\usepackage{mathrsfs}
\usepackage{amssymb}
\usepackage{amsmath}
\usepackage{amsbsy}
\usepackage{epsfig}
\usepackage{enumerate}
\usepackage{bm}
\usepackage{color}

\usepackage[colorlinks,linkcolor=blue]{hyperref}

\newtheorem{lemma}{Lemma}[section]
\newtheorem{theorem}{Theorem}[section]
\newtheorem{proposition}{Proposition}[section]

\newtheorem{definition}{Definition}[section]
\newtheorem{remark}{Remark}[section]

\def\pr{\textsf{P}} % the symbol P for probability used the sans serif letter
\def\ep{\textsf{E}} % the symbol E for expectation used the sans serif letter
\def\Sbep{\widehat{\mathbb E}} % the symbol E for sub-linear expectation
\def\cSbep{\widehat{\mathcal E}} % the symbol E for conjugate sub-linear expectation
\def\Capc{\mathbb V} % the symbol V for capacity under E
\def\upCapc{\widehat{\mathbb V}} % the symbol V with a hat  for a special upper capacity
\def\lowCapc{\widehat{\mathcal V}} % the symbol V with a hat  for a    conjugate capacity of the above
\def\outCapc{\widehat{\mathbb V}^{\ast}}% the symbol  for  the countably additive extension of the upper capacity
\def\outcCapc{\widehat{\mathcal V}^{\ast}} % the symbol  for   for conjugate countably additive capacity

\def\vSbep{\breve{\mathbb E}}

\renewcommand{\baselinestretch}{1.2}

\begin{document}

\thispagestyle{empty}

\title{\bf The moments of the maximum of normalized partial sums related to    laws of the iterated logarithm under the sub-linear expectation\thanks{ Research supported by grants from the NSF of China (Grant
No.11731012,12031005), Ten Thousands Talents Plan of Zhejiang Province (Grant No. 2018R52042), NSF of
Zhejiang Province (Grant No. LZ21A010002) and the Fundamental Research Funds for the Central Universities.}
}

\author{
\textnormal{Li-Xin Zhang}\footnote{ORCID ID:0000-0002-0121-0187;  Email:stazlx@zju.edu.cn}\\
{\sl \small School of Mathematical Sciences, Zhejiang University, Hangzhou 310027} }
\date{}
 \maketitle

\begin{abstract}
Let $\{X_n;n\ge 1\}$ be a sequence of independent and identically distributed random variables on a sub-linear expectation space $(\Omega,\mathscr{H},\Sbep)$, $S_n=X_1+\ldots+X_n$. We consider the moments of $\max_{n\ge 1}|S_n|/\sqrt{2n\log\log n}$. The sufficient and necessary
conditions for the  moments to be  finite are given. As an application, we  obtain the law of the iterated logarithm for moving average  processes of independent and identically distributed random variables.

{\bf Keywords:}  sub-linear expectation, capacity,  moments,
 laws of the iterated logarithm, moving average  process

{\bf AMS 2010 subject classifications:} 60F15; 60F05
\end{abstract}

\baselineskip 22pt

\renewcommand{\baselinestretch}{1.7}

%%%%%%%%%%%%%%%%%%%%%%%%%%%%%%%%%%
%Section 1 starts here
%%%%%%%%%%%%%%%%%%%%%%%%%%%%%%%%%%

%%%%%%%%%%%%%%%%%%%%%%%%%%%%%%%%%%%%%%%%%%%%%%%%%%%
%%%%%%%%%%%%%%%%%%%%%%%% Results %%%%%%%%%%%%%%%%%%%%

\section{ Introduction and basic settings}\label{sect1}
\setcounter{equation}{0}
      Let  $\{X_n;n\ge 1\}$ be  a sequence of independent and identically distributed (i.i.d)     random variables on a probability space $(\Omega,\mathcal F,\pr)$ with mean zeros and finite variances, $S_n=X_1+\ldots+X_n$,   $\log x=\ln \max(e,x)$ and $a_n=\sqrt{2n\log\log n}$. Siegmund (1969) and Teicher (1971) studied the moments related to  Hartman and Wintner (1941)'s law of the iterated logarithm logarithm. They obtained the sufficient and necessary conditions for  the moments of the maximum of normalized partial sums $\max_{n\ge 1}|S_n|/a_n$ to be finite. Recently, Dolera  and Regazzini  (2019) established a    reformulation of the Siegmund-Teicher inequality.
The purpose of this paper is to study the  moments  of  $\max_{n\ge 1}|S_n|/a_n$  under the sub-linear expectation.
 In the reminder of this section, we give the basic settings on the sub-linear expectations.  The main results will be stated in section \ref{sectMain}, and the upper bound of  the Siegmund-Teicher type inequality established by   Dolera  and Regazzini  (2019) is improved to the optimal one.  In section \ref{sectAppl},  as an application of the main results, the law of iterated logarithm for moving average processes generated by a sequence of   i.i.d. random variables is established. The proof of the main results is given in the last section.

We use the framework and notations of Peng (2008,2019).   Let  $(\Omega,\mathcal F)$
 be a given measurable space  and let $\mathscr{H}$ be a linear space of real functions
defined on $(\Omega,\mathcal F)$ such that if $X_1,\ldots, X_n \in \mathscr{H}$,  then $\varphi(X_1,\ldots,X_n)\in \mathscr{H}$ for each
$\varphi\in C_{l,Lip}(\mathbb R_n)$,  where $C_{l,Lip}(\mathbb R_n)$ denotes the linear space of (local Lipschitz)
functions $\varphi$ satisfying
\begin{eqnarray*} & |\varphi(\bm x) - \varphi(\bm y)| \le  C(1 + |\bm x|^m + |\bm y|^m)|\bm x- \bm y|, \;\; \forall \bm x, \bm y \in \mathbb R_n,&\\
& \text {for some }  C > 0, m \in \mathbb  N \text{ depending on } \varphi. &
\end{eqnarray*}
 We also denote $C_{b,Lip}(\mathbb R^n)$ the space of bounded  Lipschitz
functions.

\begin{definition}\label{def1.1} A  sub-linear expectation $\Sbep$ on $\mathscr{H}$  is a function $\Sbep: \mathscr{H}\to \overline{\mathbb R}$ satisfying the following properties: for all $X, Y \in \mathscr H$, we have
\begin{description}
  \item[\rm (a)]   if $X \ge  Y$ then $\Sbep [X]\ge \Sbep [Y]$;
\item[\rm (b)]  $\Sbep[X+Y]\le \Sbep [X] +\Sbep [Y ]$ whenever $\Sbep [X] +\Sbep [Y ]$ is not of the form $+\infty-\infty$ or $-\infty+\infty$;
\item[\rm (c)] $\Sbep [c] = c$ and $\Sbep [\lambda X] = \lambda \Sbep  [X]$, $\lambda>0$.
 \end{description}
 Here $\overline{\mathbb R}=[-\infty, \infty]$. The triple $(\Omega, \mathscr{H}, \Sbep)$ is called a sub-linear expectation space. Given a sub-linear expectation $\Sbep $, let us denote the conjugate expectation $\cSbep$of $\Sbep$ by
$ \cSbep[X]:=-\Sbep[-X]$, $\forall X\in \mathscr{H}$.
\end{definition}

\begin{definition} ({\em Peng (2008, 2019)})
\begin{description}
  \item[ \rm (i)]  Let $\bm X_1$ and $\bm X_2$ be two $n$-dimensional random vectors defined,
respectively, on sub-linear expectation spaces $(\Omega_1, \mathscr{H}_1, \Sbep_1)$
  and $(\Omega_2, \mathscr{H}_2, \Sbep_2)$. $\bm X_1$ and $\bm X_2$ are said  to be  identically distributed,  denoted by $\bm X_1\overset{d}= \bm X_2$,  if
$$ \Sbep_1[\varphi(\bm X_1)]=\Sbep_2[\varphi(\bm X_2)], \;\; \forall \varphi\in C_{b,Lip}(\mathbb R_n). $$
 A sequence $\{X_n;n\ge 1\}$ of random variables is said to be identically distributed if $X_i\overset{d}= X_1$ for each $i\ge 1$.
\item[\rm (ii)]     On a sub-linear expectation space  $(\Omega, \mathscr{H}, \Sbep)$, a random vector $\bm Y =
(Y_1, \ldots, Y_n)$, $Y_i \in \mathscr{H}$ is said to be independent relative to another random vector $\bm X =
(X_1, \ldots, X_m)$ , $X_i \in \mathscr{H}$ under $\Sbep$ if
$ \Sbep [\varphi(\bm X, \bm Y )] = \Sbep \big[\Sbep[\varphi(\bm x, \bm Y )]\big|_{\bm x=\bm X}\big]$   for each  $\varphi\in C_{b,Lip}(\mathbb R_m \times \mathbb R_n)$.
 A sequence  $\{X_n; n\ge 1\}$ of random variables
 is said to be independent, if
 $X_{i+1}$ is independent relative to $(X_{i+1},\ldots, X_n)$ for each $i\ge 1$.
 \end{description}
\end{definition}

 Let $(\Omega, \mathscr{H}, \Sbep)$ be a sub-linear space.  We denote a pair $(\upCapc,\lowCapc)$ of capacities by
\begin{equation}\label{equpCapc} \upCapc(A):=\inf\{\Sbep[\xi]: I_A\le \xi, \xi\in\mathscr{H}\}, \;\; \lowCapc(A):= 1-\upCapc(A^c),\;\; \forall A\in \mathcal F,
\end{equation}
where $A^c$  is the complement set of $A$. It is obvious that $\upCapc$ is a sub-additive capacity in sense that $\upCapc(A\cup B)\le \upCapc(A)+\upCapc(B)$, and
\begin{equation}\label{eq1.3}
 \Sbep[f]\le \upCapc(A)\le \Sbep[g], \;\;\cSbep[f]\le \lowCapc(A) \le \cSbep[g],\;\;
\text{ if } f\le I_A\le g, f,g \in \mathscr{H}.
\end{equation}

Also, we define the  Choquet integrals/expectations $(C_{\upCapc},C_{\lowCapc})$  by
$$ C_V[X]=\int_0^{\infty} V(X\ge t)dt +\int_{-\infty}^0\left[V(X\ge t)-1\right]dt, $$
with $V$ being replaced by $\upCapc$ and $\lowCapc$, respectively.

Because the capacity $\upCapc$ defined as in \eqref{equpCapc}  may be not countably sub-additive,   we  consider its countably sub-additive extension.
\begin{definition}   A  countably sub-additive extension $\outCapc$  of $\upCapc$   is defined by
\begin{equation}\label{outcapc} \outCapc(A)=\inf\Big\{\sum_{n=1}^{\infty}\upCapc(A_n): A\subset \bigcup_{n=1}^{\infty}A_n, \; A_n\in \mathcal F\Big\},\;\; \outcCapc(A)=1-\outCapc(A^c),\;\;\; A\in\mathcal F.
\end{equation}
\end{definition}

Finally, for real numbers $x$ and $y$, denote $x\vee y=\max(x,y)$, $x\wedge y=\min(x,y)$, $x^+=\max(0,x)$ and $x^-=\max(0,-x)$.
For a random variable $X$, because $XI\{|X|\le c\}$  may be not in $\mathscr{H}$, we will truncate $X$ in the form $(-c)\vee X\wedge c$ denoted by $X^{(c)}$.
 Denote $\vSbep[X]=\lim_{c\to \infty} \Sbep[X^{(c)}]$  if the limit exists. It can be verified that $\vSbep[|X|]\le C_{\upCapc}(|X|)$.

\section{Main results}\label{sectMain}
Let $(\Omega, \mathscr{H}, \Sbep)$ be a sub-linear expectation space and    $\{X, X_n;n\ge 1\}$ be a sequence  random variables on it.
 Denote $S_n=\sum_{i=1}^n X_i$,   $a_n=\sqrt{2n\log\log n}$.

When $X_1,X_2,\ldots$ are i.i.d. random variables on a classical probability space $(\Omega,\mathcal F, \pr)$, Siegmund (1969) and Teicher (1971) studied the moments of $\max_{n\ge 1}|S_n|/a_n$. Particularly, it is shown that $\ep\left[\max_{n\ge 1}|S_n|/a_n\right]<\infty$ if and only if
$\ep[X]=0$ and $\ep[X^2]<\infty$ (see (16) of Siegmund (1969)). Dolera  and Regazzini  (2019) established a    reformulation of the Siegmund-Teicher inequality for $\ep\left[\max_{n\ge 1}|S_n|^r/a_n^r\right]$ ($r>2)$.
In this paper, we consider the same problem under the sub-linear expectation. Since $\max_{n\ge 1}|S_n|/a_n$ may be not an element of $\mathscr{H}$ and so its sub-linear expectation has no definition, and $\Sbep[|S_n|/a_n]$ may be not finite under the conditions for the law of iterated logarithm, we consider the  Choquet  expectation  of  $\max_{1\le n\le N} |S_n|/a_n$. The following is our main result:
   \begin{theorem} \label{thMomentLIL}  Let $\{X_n; n\ge 1\}$  be a sequence of independent random variables   with $X_n\overset{d}= \alpha_nX$, $|\alpha_n|\le 1$. Denote $\overline{\sigma}_X^2=\breve{\mathbb E}[X^2]$,
   \begin{equation}\label{eqthMLIL.1} \varsigma_X\hat{=}C_{\upCapc}\left[\frac{X^2}{\log\log|X|}\right],
    \end{equation}
     \begin{equation}\label{eqthMLIL.0} \eta_{X,r}\hat{=}\begin{cases} \varsigma_X, &\text{ if } 0<r<2,\\
     C_{\upCapc}\left[\frac{X^2\log |X|}{\log\log|X|}\right], &\text{ if } r=2,\\
     C_{\upCapc}\left[|X|^r\right], & \text{ if } r>2.
     \end{cases}
    \end{equation}
  \begin{description}
   \item[\rm (a)] Suppose
   $\breve{\mathbb E}[X_n]\le 0$.
       Then, for any $r>0$ and $p>2\vee r$, there exists a constant $K_{r,p}$ such that
\begin{equation}\label{eqthMLIL.3}
   \max_{N}C_{\upCapc}\left[\max_{1\le n\le N}\left(\frac{ S_n^+ }{a_n}\right)^r \right]
 \le  K_{r,p}  \left\{ \eta_{X,r}+\varsigma_X^{r/p}+\overline{\sigma}_X^r\right\};
\end{equation}
and, for any $r>2$, there exists a constant $K_r$ such that
\begin{equation}\label{eqthMLIL.4}
   \max_{N}C_{\upCapc}\left[\max_{1\le n\le N}\left(\frac{ S_n^+ }{a_n}\right)^r \right]
 \le  K_r C_{\upCapc}[|X|^r].
\end{equation}
\item[\rm (b)] Suppose   $X_n\overset{d}= X$, $r>0$, and
   \begin{equation}\label{eqthMLIL.5}
 \max_{N}C_{\upCapc}\left[\max_{1\le n\le N}\frac{ |S_n|^r }{a_n^r}\right]<\infty.
\end{equation}
Then,   $\overline{\sigma}_X^2<\infty$,
   $\breve{\mathbb E}[X]=\breve{\mathbb E}[-X]=0$ and $\eta_{X,r}<\infty$.
\end{description}
\end{theorem}

It is obvious that when $r\ge 2$, $\eta_{X,r}<\infty$ implies that $\varsigma_X<\infty$ and $\overline{\sigma}_X<\infty$, so the results coincide with those of Siegmund (1969) and Teicher (1971). When $r>2$, \eqref{eqthMLIL.4} has improved Proposition 2.1 of Dolera  and Regazzini  (2019). In this case, the bound in the right hand of \eqref{eqthMLIL.4} is optimal, because the left hand of \eqref{eqthMLIL.5} does not exceed $C_{\upCapc}[|X|^r]/a_1^r$.    When $0<r<2$, the sufficient and necessary conditions for \eqref{eqthMLIL.5}  coincide with those for the law of the iterated logarithm (see Zhang (2021)). It is interesting that if  \eqref{eqthMLIL.5} holds for one $0<r<2$, then it holds for all $0<r<2$.

%%%%%%%%%%%%%%%%%%%%%%%%%%%%%%%%%%%%%%%%%%%%%%%%%%
  \section{Applications}\label{sectAppl}
  \setcounter{equation}{0}
 As an application of Theorem \ref{thMomentLIL}, we show the law of the iterated logarithm for the moving average process in this section.
Let $\{Y_i;i\ge 1\}$ be a sequence of i.i.d. random variables under the sub-linear expectation $\Sbep$. We consider the moving average process
\begin{equation}\label{eqmoving1}X_t=\sum_{j=-\infty}^{t-1} \beta_j Y_{t-j},
\end{equation}
where
\begin{equation}\label{eqmoving2} B\hat{=}\sum_{j=-\infty}^{\infty}|\beta_j|<\infty, \;\; \beta\hat{=}\sum_{j=-\infty}^{\infty}\beta_j.
\end{equation}
Let $T_n=\sum_{t=1}^n X_t$. Define $Y_t=0$ for $t=0,-1,-2,\ldots$. Then
\begin{equation}\label{eqmoving3} X_t= \sum_{j=-\infty}^{\infty} \beta_j Y_{t-j}, \;\; T_n =\sum_{j=-\infty}^{\infty} \beta_j\sum_{t=1}^n Y_{t-j}.
\end{equation}
If  $\{Y_i\}_{i=-\infty}^{\infty}$ is a bi-directional sequence of i.i.d. random variables, i.e., for any $i_1<i_2<\cdots<i_p$, $\{Y_{i_1},\ldots,Y_{i_p}\}$ are i.i.d. random variables, we can define $X_t$ as in \eqref{eqmoving3} directly.  To see that $X_t$ is well-defined, we need to show that the infinite series in \eqref{eqmoving1} or \eqref{eqmoving3} are convergent. This is not a problem when $\Sbep=\ep$ is a linear mathematical expectation and $\ep[|Y_1|]<\infty$, since $\ep[|X_t|]\le \sum_{j=-\infty}^{\infty}|\beta_j| \ep[|Y_{t-j}|]<\infty$. However, in the sub-linear expectation space, first, the sub-linear expectation $\Sbep$ may be not countably sub-additive, and secondly, though every $Y_j\in \mathscr{H}$,  $X_t$ may be not an element of $\mathscr{H}$   so that $\Sbep[|X_t|]$ may have no definition.
 To consider the convergence of a random infinite series, we need a lemma.
\begin{lemma}\label{lem5.1} Let $\{\xi_j;j=0,\pm 1, \pm 2,\ldots\}$ be a sequence of random variables on $(\Omega,\mathscr{H},\Sbep)$ with
$ \sum_{j=-\infty}^{\infty}\vSbep[|\xi_j|]<\infty.
$ Then,
\begin{equation}\label{eqlem5.1.2} \outcCapc\big(\sum_{j=-\infty}^{\infty}|\xi_j|<\infty\big)=1 \end{equation}
and
\begin{equation}\label{eqlem5.1.3} \outCapc\big(\sum_{j=-\infty}^{\infty}|\xi_j|\ge x\big)\le x^{-1}\sum_{j=-\infty}^{\infty}\vSbep[|\xi_j|], \;\; x>0. \end{equation}
\eqref{eqlem5.1.3} is obvious when $\sum_{j=-\infty}^{\infty}\vSbep[|\xi_j|]=\infty$.
\end{lemma}
{\bf Proof.} Let $M_n=\sum_{|j|\le n}|\xi_j|$. Then,
\begin{align*}
& \upCapc(M_{n+m}-M_n\ge \epsilon)=\upCapc\left(\sum_{|j|=n+1}^{n+m}(|\xi_j|\wedge \epsilon)\ge \epsilon\right) \\
\le & \frac{1}{\epsilon}\sum_{|j|=n+1}^{n+m}\Sbep[|\xi_j|\wedge\epsilon]\le \frac{1}{\epsilon}\sum_{|j|=n+1}^{n+m}\vSbep[|\xi_j|]\to 0\; \text{ as } n\to \infty.
\end{align*}
For $\epsilon_k=2^{-k}$, there exists a sequence $n_k\nearrow \infty$ such that
\begin{equation}\label{eqlem5.1.4} \upCapc(M_{n_{k+1}}-M_{n_k}\ge \epsilon_k)\le \epsilon_k.
\end{equation}
It follows that
$ \sum_{k=1}^{\infty}\outCapc(M_{n_{k+1}}-M_{n_k}\ge \epsilon_k)<\infty.
$
Hence $\outCapc(A)=0$, where $A=\{M_{n_{k+1}}-M_{n_k}\ge \epsilon_k\; i.o.)$.
On $A^c$,
$$ \sum_{j=-\infty}^{\infty}|\xi_j|=M_{n_1}+\sum_{k=1}^{\infty}(M_{n_{k+1}}-M_{n_k})<M_{n_1}
+\sum_{k=1}^{\infty}\epsilon_k<\infty. $$
\eqref{eqlem5.1.2} is proven. Furthermore, by \eqref{eqlem5.1.4},
\begin{align*}
& \outCapc\left(\sum_{|j|\ge n_k+1}|\xi_j|\ge 2\epsilon_k\right)
=\outCapc\left(\sum_{l=k}^{\infty} (M_{n_{l+1}}-M_{n_l})\ge \sum_{l=k}^{\infty}\epsilon_l\right) \\
&\;\; \le   \sum_{l=k}^{\infty} \outCapc\left( M_{n_{l+1}}-M_{n_l}\ge  \epsilon_l\right)
<\sum_{l=k}^{\infty} \epsilon_l=2\epsilon_k.
\end{align*}
For $x>0$, choose $k$ such that $2\epsilon_k<x$. Then by noting that $I\{M_{n_k}\ge y\}\le y^{-1}(|M_{n_k}|\wedge y)$ and
$|M_{n_k}|\wedge y\in \mathscr{H}$, we have that
\begin{align*}
& \outCapc\left(\sum_{j=-\infty}^{\infty}|\xi_j|\ge x\right)
\le   \outCapc\left(\sum_{|j|\ge n_k+1}|\xi_j|\ge 2\epsilon_k\right)+\outCapc\left(M_{n_k}\ge x-2\epsilon_k\right)\\
&\;\; \le  2\epsilon_k+\frac{1}{x-2\epsilon_k}\vSbep[|M_{n_k}|]
\le 2\epsilon_k+\frac{1}{x-2\epsilon_k}\sum_{j=-\infty}^{\infty} \vSbep[|\xi_j|].
\end{align*}
Letting $k\to \infty$ yields \eqref{eqlem5.1.3}. \hfill $\Box$

By Lemma \ref{lem5.1}, the infinite series on the right hands of \eqref{eqmoving1} and \eqref{eqmoving3} are convergent  when $\vSbep[|Y_1|]<\infty$,  so $X_t$ and $T_n$ are well-defined  though they may be not   elements of $\mathscr{H}$.   The following theorems give  the law of the iterated logarithm for $\{X_t; t\ge 1\}$:

\begin{theorem} \label{thLILMA} Let $\{Y_i;i\ge 1\}$ be a sequence of  i.i.d. random variables under the sub-linear expectation $\Sbep$. Suppose that
\begin{equation}\label{eqthLILMA.1}
\breve{\mathbb E}[Y_1]=\breve{\mathbb E}[-Y_1]=0,
\end{equation}
\begin{equation}\label{eqthLILMA.2}\overline{\sigma}_Y^2=\vSbep[Y_1^2]<\infty
\end{equation} and
\begin{equation}\label{eqthLILMA.3}\varsigma_Y\hat{=}C_{\upCapc}\left[\frac{Y_1^2}{\log\log|Y_1|}\right]<\infty.
\end{equation}
 Then
 \begin{equation}\label{eqthLILMA.4}\outcCapc\left(\limsup_{n\to \infty}\frac{|T_n|}{a_n}\le |\beta| \overline{\sigma}_Y\right)=1.
 \end{equation}
 Furthermore, suppose that $\Sbep$ is regular on $\mathscr{H}_b=\{\varphi\in \mathscr{H};\varphi\text{ is bounded }\}$ in the sense that $\Sbep[\varphi_n]\to 0$ whenever $\mathscr{H}_b\ni\varphi\searrow 0$. Then
 $$\outCapc\left(C\Big\{\frac{T_n}{a_n}\Big\}=\big[-|\beta|\overline{\sigma}_Y,|\beta|\overline{\sigma}_Y\big]\right)=1. $$
 \end{theorem}
\begin{remark}
Recently, Liu and Zhang proved \eqref{eqthLILMA.4} under a condition that $\Sbep[Y_1^2(\log |Y_1|)^{1+\delta}]<\infty$ for some $\delta>0$ which is much more stringent than \eqref{eqthLILMA.2} and \eqref{eqthLILMA.3}. When $\beta_0=1$ and $\beta_j=0$ for all $j\ne 0$, $X_t=Y_t$. At this case, by Theorem 5.9 of Zhang (2021), \eqref{eqthLILMA.1}-\eqref{eqthLILMA.3} are also necessary conditions for the law of the iterated  logarithm.
\end{remark}

Theorems \ref{thLILMA}    follows  directly from Theorem 5.9 of Zhang (2021) on the law of the iterated logarithm for i.i.d. random variables,      by noting the following approximation.

\begin{proposition}Let $\{Y_i;i\ge 1\}$ be a sequence of i.i.d. random variables under the sub-linear expectation $\Sbep$ satisfying \eqref{eqthLILMA.1}-\eqref{eqthLILMA.3}.
Then
\begin{equation}\label{eqprop1}
\outcCapc\left(\lim_{n\to \infty}\frac{T_n-\beta\sum_{t=1}^n Y_t}{a_n}=0\right)=1.
\end{equation}
\end{proposition}
{\bf Proof.} Let
$ T_{m,n}= \sum_{t=1}^n \sum_{j=-m}^m \beta_j Y_{t-j}. $
Then
$T_n=T_{m,n}+ \sum_{t=1}^n \sum_{|j|>m} \beta_j Y_{t-j}. $
By Theorem \ref{thMomentLIL} (with $r=1$,$p=3$),
$$ \vSbep\left[\max_{N_1\le n\le N_2}\frac{|\sum_{t=1}^n Y_{t-j}|}{a_n}\right]\le K \left\{\varsigma_Y^{1/3}+\varsigma_Y+\overline{\sigma}_Y\right\}:=C_0. $$
By Lemma \ref{lem5.1}, it follows that
\begin{align*}
&\outCapc\left(\max_{N_1\le n\le N_2}\frac{|T_n-T_{m,n}|}{a_n}\ge \epsilon\right)
\le   \outCapc\left(\sum_{|j|>m} |\beta_j|\max_{N_1\le n\le N_2}\frac{|\sum_{t=1}^n Y_{t-j}|}{a_n}\ge \epsilon\right)\\
\le & \epsilon^{-1}\sum_{|j|>m} |\beta_j|\vSbep\left[\max_{N_1\le n\le N_2}\frac{|\sum_{t=1}^n Y_{t-j}|}{a_n}\right]
\le C_0 \epsilon^{-1}\sum_{|j|>m} |\beta_j|.
\end{align*}
Notice that
\begin{align*}
 T_{m,n}= \sum_{j=-m}^m \beta_j \sum_{t=1}^{n} Y_{t}+\sum_{j=1}^m \beta_j\sum_{t=1-j}^0(Y_t-Y_{n+t})
+\sum_{j=1}^m \beta_{-j}\sum_{t=1}^j(Y_t-Y_{n+t}).
\end{align*}
It follows that
\begin{align*}
|T_{m,n}-\beta\sum_{t=1}^n Y_t|\le
\sum_{|j|>m}|\beta_j| |\sum_{t=1}^n Y_t|+B\sum_{t=-m}^m |Y_t|+B\sum_{t=-m}^m |Y_{n+t}|.
\end{align*}
By Lemma \ref{lem5.1} and Theorem \ref{thMomentLIL}, we have that
\begin{align*}
&\outCapc\left(\sum_{|j|>m}|\beta_j| \max_{N_1\le n\le N_2} \frac{|\sum_{t=1}^n Y_t|}{a_n}\ge \epsilon\right) \\
 \le   &\epsilon^{-1}\sum_{|j|>m}|\beta_j|\vSbep\left[\max_{N_1\le n\le N_2} \frac{|\sum_{t=1}^n Y_t|}{a_n}\right]
\le \epsilon^{-1} C_0 \sum_{|j|>m}|\beta_j|, \\
& \outCapc\left(B\sum_{t=-m}^m |Y_t|\ge \epsilon a_n\right)\le \frac{B\sum_{t=-m}^m\vSbep[|Y_t|]}{\epsilon a_n}
=\frac{2mB\vSbep[|Y_1|]}{\epsilon a_n }, \\
& \outCapc\left(B\sum_{t=-m}^m \max_{N_1\le n\le N_2}\frac{|Y_{n+t}|}{a_n}\ge \epsilon\right)\\
 \le & \sum_{t=-m}^m \sum_{  N_1\le n\le N_2}\outCapc\left(B  \frac{|Y_{n+t}|}{a_n}\ge \epsilon/m\right)
\le 2m\sum_{  n\ge N_1}\upCapc\left(|Y_1|\ge \epsilon a_n/(2Bm)\right).
\end{align*}
Combing the above arguments yields that
\begin{align*}
&\max_{N_2}\outCapc\left(\max_{N_1\le n\le N_2}\frac{|T_n-\beta\sum_{t=1}^n Y_t|}{a_n}\ge 4\epsilon\right) \\
\le & \frac{2C_0\sum_{|j>m}|\beta_j|}{\epsilon}+\frac{2mB\vSbep[|Y_1|]}{\epsilon a_{N_1} }
+2m\sum_{  n\ge N_1}\upCapc\left(|Y_1|\ge \epsilon a_n/(Bm)\right) \\
& \; \to 0\;\; \text{ as } N_1\to 0 \text{ and then } m\to \infty.
\end{align*}
Now, let $\epsilon_k=2^{-k}$. Then there exists a sequence $n_k\nearrow \infty$ such that
$$ \outCapc\left(\max_{n_k\le n\le n_{k+1}}\frac{|T_n-\beta\sum_{t=1}^n Y_t|}{a_n}\ge \epsilon_k\right)\le \epsilon_k. $$
It follows that $\sum_{k=1}^{\infty}\outCapc(A_k)<\infty$, where $A_k=\left\{\max_{n_k\le n\le n_{k+1}}\frac{|T_n-\beta\sum_{t=1}^n Y_t|}{a_n}\ge \epsilon_k\right\}$. Hence, by the countable sub-additivity of $\outCapc$, $\outCapc\left(A_k\;\; i.o.\right)=0$. On the event  $\left\{A_k\;\; i.o.\right\}^c$,
$$ \limsup_{n\to\infty}\frac{|T_n-\beta\sum_{t=1}^n Y_t|}{a_n}\le
\limsup_{k\to\infty}\max_{n_k\le n\le n_{k+1}}\frac{|T_n-\beta\sum_{t=1}^n Y_t|}{a_n}=0. $$
\eqref{eqprop1} is proven. \hfill $\Box$.

%%%%%%%%%%%%%%%%%%%%%%%%%%%%%%%%%%%%%%%%%%%%%%%%%%
  \section{Proofs}\label{sectProof}
  \setcounter{equation}{0}

   For proving the main results, we need several inequalities.
\begin{lemma}\label{lemExpIneq}    Suppose that  $\{X_1,\ldots, X_n\}$ is a sequence  of independent random variables
on $(\Omega, \mathscr{H}, \Sbep)$. Set $S_n=\sum_{k=1}^n X_k$,
$A_n(p,y)=\sum_{i=1}^n\Sbep[(X_i^+\wedge y )^p]$ and
$ \breve{B}_{n,y}=\sum_{i=1}^n \vSbep[(X_i\wedge y)^2]$.    Then,
   for all $p\ge 2$, $x,y>0$, $0<\delta\le 1$,
\begin{align}\label{eqExpIneq.5}
& \upCapc\Big( \max_{k\le n}  \sum_{i=1}^k(X_i-\vSbep[X_i])\ge x\Big) \nonumber\\
\le &  \upCapc\big(\max_{k\le n} X_k> y \big)
  +2\exp\{p^p\}\Big\{\frac{A_n(p,y)}{y^p} \Big\}^{\frac{\delta x}{10y}}
+\exp\left\{-\frac{x^2}{2\breve{B}_{n,y}(1+\delta)   }\right\}.
\end{align}
  \end{lemma}

{\bf Proof}.   The proof of \eqref{eqExpIneq.5} is the same as that of (3.1) of  Zhang (2021) if we note
$$ \Sbep[e^{t (X_k\wedge y)}]=\vSbep[e^{t (X_k\wedge y)}]\le 1+t\vSbep[X_k]+\frac{e^{ty}-1-t y}{y^2} \vSbep[(X_k\wedge y)^2], \; y>0. \;\;\;\; \Box$$

\begin{lemma} \label{lem3} Suppose $X\in \mathscr{H}$, $r>0$. Let  $\varsigma_X$ and $\eta_{X,r}$ be defined as in \eqref{eqthMLIL.1} and \eqref{eqthMLIL.0}, respectively. Then,
\begin{description}
  \item[\rm (i)]  for any $\delta>0$ and $p>2$,
$$ \sum_{n=1}^{\infty} \frac{C_{\upCapc}\big[\big(|X|\wedge (\delta a_n)\big)^p\big]}{a_n^p}\le c_{\delta,p}\varsigma_X; $$
 \item[\rm (ii)] for any $r>0$,
$$ \sum_{n=1}^{\infty} \frac{C_{\upCapc}\big[\big((|X|-a_n\big)^+\big)^r\big]}{a_n^r}\le c_r \eta_{X,r}. $$
\end{description}
\end{lemma}
{\bf Proof}. (i) Let $ f(x)$ be the inverse function of $\sqrt{2x\log\log x}$. Then,
\begin{align*}
  C_{\upCapc} \big[\big(|X|\wedge (\delta a_n)\big)^p\big]\le  &\int_0^{(\delta a_n)^p}\upCapc\left(|X|^p>x\right)dx
=p\int_0^{\delta a_n}x^{p-1}\upCapc\left(|X|>x\right)dx \\
\le & 2p\int_0^{(\delta^2+1) n} (2\log\log y)^{\frac{p}{2}} y^{\frac{p}{2}-1} \upCapc\left(f(|X|)>y\right)dy.
\end{align*}
It follows that
\begin{align*}
&\sum_{n=1}^{\infty} \frac{C_{\upCapc}\big[\big(|X|\wedge (\delta a_n)\big)^p\big]}{a_n^p}
\le  2p\sum_{n=1}^{\infty}  a_n^{-p}\int_0^{(\delta^2+1) n} (2\log\log y)^{\frac{p}{2}} y^{\frac{p}{2}-1} \upCapc\left(f(|X|)>y\right)dy\\
\le &  22^pp\int_{0}^{\infty}(2x\log\log x)^{-p/2} \int_0^{(\delta^2+1) x} (2\log\log y)^{\frac{p}{2}} y^{\frac{p}{2}-1} \upCapc\left(f(|X|)>y\right)dydx \\
\le &  22^pp\int_{0}^{\infty} (2\log\log y)^{\frac{p}{2}} y^{\frac{p}{2}-1} \upCapc\left(f(|X|)>y\right)dy\int_{y/(\delta^2+1)}^{\infty}  (2x\log\log x)^{-p/2}dx \\
\le &c_{\delta,p}\int_{0}^{\infty}  \upCapc\left(f(|X|)>y\right)dy=c_{\delta,p}C_{\upCapc}\big[f(|X|)\big]
\le c_{\delta,p}C_{\upCapc}\left[\frac{|X|^2}{\log\log |X|}\right].
\end{align*}
For (ii), we have that
\begin{align*}
&  \sum_{n=17}^{\infty} \frac{C_{\upCapc}\big[\big(\big(|X|-a_n\big)^+\big)^r\big]}{a_n^r}
\le   \sum_{n=17}^{\infty} \frac{\int_{a_n}^{\infty}rx^{r-1}\upCapc(|X|>x)dx}{a_n^r}\\
\le & \int_{16}^{\infty} \frac{\int_{a_y}^{\infty}rx^{r-1}\upCapc(|X|>x)dx}{a_y^r} dy
=\int_{a_{16}}^{\infty} rx^{r-1}\upCapc(|X|>x)dx  \int_{16\le y, a_y\le x} \frac{1}{a_y^r} dy.
\end{align*}
When $0<r<2$,
$$ \int_{16\le y, a_y\le x} \frac{1}{a_y^r} dy\approx   \frac{1}{x^r}\frac{x^2}{\log\log x}, $$
so
$$ \sum_{n=17}^{\infty} \frac{C_{\upCapc}\big[\big(\big(|X|-a_n\big)^+\big)^r\big]}{a_n^r}\le
c_r\int_{a_{16}}^{\infty} \frac{x }{\log\log x} \upCapc(|X|>x)dx \le c_rC_{\upCapc}\left[\frac{X^2}{\log\log |X|}\right]. $$
When $r=2$,
$$ \int_{16\le y, a_y\le x} \frac{1}{a_y^r} dy\approx  \frac{\log x}{\log\log x}, $$
 so
$$ \sum_{n=17}^{\infty} \frac{C_{\upCapc}\big[\big(\big(|X|-a_n\big)^+\big)^r\big]}{a_n^r}\le
c_r\int_{a_{16}}^{\infty} \frac{x\log x }{\log\log x} \upCapc(|X|>x)dx \le c_rC_{\upCapc}\left[\frac{X^2\log |X|}{\log\log |X|}\right]. $$
When $r>2$,
$$ \int_{16\le y, a_y\le x} \frac{1}{a_y^r} dy\approx Const., $$
 so
$$ \sum_{n=17}^{\infty} \frac{C_{\upCapc}\big[\big(\big(|X|-a_n\big)^+\big)^r\big]}{a_n^r}\le
c_r\int_{a_{16}}^{\infty} rx^{r-1}\upCapc(|X|>x)dx \le c_rC_{\upCapc}\left[|X|^r\right]. $$
Furthermore,
$$ \sum_{n=1}^{16} \frac{C_{\upCapc}\big[\big((|X|-a_n)^+\big)^r\big]}{a_n}
\le 16 C_{\upCapc}\big[\big((|X|-\sqrt{2})^+\big)^r\big]\le c_r \eta_{X,r}. $$
The proof is complete, and, with a similar argument it can be shown that
\begin{equation}\label{eqprooflem3.5} \int_{16}^{\infty} \frac{\int_{a_y}^{\infty}rx^{r-1}\upCapc(|X|>x)dx}{a_y^r} dy<\infty \implies \eta_{X,r}<\infty. \;\;\;\; \Box
\end{equation}

\begin{lemma}\label{lem4} Suppose that for some $C_0>0$, $p>0$,
$x^p\upCapc(|X|\ge x)\le C_0$, $x>0$.
Then, for $0<r<p$,
$$C_{\upCapc}[|X|^r]=\int_0^{\infty} \upCapc(|X|\ge x^{1/r})dx \le \int_0^{\infty} 1\wedge (C_0x^{-p/r}) dx= \frac{r}{p-r}C_0^{r/p}. $$
\end{lemma}
{\bf Proof.} Trivial. \hfill $\Box$

{\bf Proof of Theorem \ref{thMomentLIL}.} Write $\overline{\sigma}=\overline{\sigma}_X$, $\varsigma=\varsigma_X$, $\eta_r=\eta_{X,r}$, and suppose that they are finite.
Let $a_y=\sqrt{2y\log\log y}$.

(a)
  Let $n_k=2^k$, $I(k)=\{n_k+1,\ldots, n_{k+1}\}$.  Then
 \begin{align*}
   \max_{0\le k\le m}\max_{n_k\le  n\le n_{k+1}} \frac{S_n} {a_n}
\le   2\max_{0\le k\le m} \frac{\max\limits_{1\le  n\le n_{k+1}}S_n^+} {a_{n_{k+1}}}=:2\max_{0\le k\le m}I_k.
\end{align*}
Applying \eqref{eqExpIneq.5} with $\delta=1$, $p>2\vee r$, $x=:x_{k+1}=z^{1/r}a_{n_{k+1}}$ and $y=y_{k+1}=x/30$ yields
\begin{align}\label{eq:ad}
\upCapc(&\max_{0\le k\le m}I_k\ge z^{1/r})\le   \sum_k \upCapc\left(\max_{n\le n_{k+1}}X_n\ge z^{1/r} a_{n_{k+1}}/30\right)\nonumber \\
&+c_p \sum_k \left(\frac{n_{k+1}\vSbep[(|X|\wedge y_{k+1})^p]}{y_{k+1}^p} \right)^3
 +\sum_k\exp\left\{-\frac{z^{2/r}\log\log n_{k+1}}{4\overline{\sigma}^2}\right\}\nonumber \\
 &   =:g_1(z)+g_2(z)+g_3(z).
\end{align}
Note
\begin{align*}
 g_1(z)\le &\sum_kn_{k+1}\upCapc\left(X\ge z^{1/r} a_{n_{k+1}}/30\right)
\le   2\sum_{n=1}^{\infty} \upCapc\left(|X|\ge z^{1/r} a_n/30\right) \\
\le &2 \sum_{n=1}^{\infty} \left(\upCapc\left(|X|\wedge  a_n\ge z^{1/r} a_n/60\right)
+     \upCapc\left((|X|- a_n)^+\ge z^{1/r} a_n/60\right)\right) \\
\le & K_{r,p} z^{-p/r}\sum_{n=1}^{\infty}\frac{\Sbep[(|X|\wedge a_n)^p]}{a_n^p}+2\sum_{n=1}^{\infty}\upCapc\left((|X|- a_n)^+\ge z^{1/r}a_n/60\right).
\end{align*}
By Lemmas \ref{lem4} and \ref{lem3}, we have that
\begin{align*}
\int_0^{\infty} 1\wedge g_1(z) dz \le &  K_{r,p} \left(\sum_{n=1}^{\infty}\frac{\Sbep[(|X|\wedge a_n)^p]}{a_n^p}\right)^{r/p}+K_r \sum_{n=1}^{\infty} \frac{C_{\upCapc}[((|X|- a_n)^+)^r]}{a_n^r}\\
\le & K_{r,p} \varsigma^{r/p}+K_r\eta_r.
\end{align*}
  For $g_2(z)$, note $ \Sbep[(|X|\wedge y_{k+1})^p]/y_{k+1}^p\le c_p\Sbep[(|X|\wedge a_{n_{k+1}})^p]/y_{k+1}^p+c_{p}
C_{\upCapc}[((|X|-a_{n_{k+1}})^+)^r]/y_{k+1}^r$. By Lemma  \ref{lem3}, it follows that
\begin{align*}
g_2(z)\le & K_{r,p}  \sum_k \left\{ \Big(\sum_{n\in I(k)} \frac{\Sbep[(|X|\wedge a_n)^p]}{z^{p/r}a_n^p}\Big)^3+     \Big(\sum_{n\in I(k)} \frac{C_{\upCapc}[((|X|-a_n)^+)^r]}{z a_n^r}\Big)^3\right\}\\
\le & K_{r,p}z^{-3p/r}\Big(\sum_{n=1}^{\infty} \frac{\Sbep[(|X|\wedge a_n)^p]}{a_n^p}\Big)^3+K_{r,p}z^{-3}   \Big(\sum_{n=1}^{\infty} \frac{C_{\upCapc}[((|X|-a_n)^+)^r]}{a_n^r}\Big)^3\\
\le & K_{r,p}z^{-3p/r}\varsigma^3+K_{r,p}z^{-3}\eta_r^3.
\end{align*}
Hence, by Lemma \ref{lem4},
$ \int_0^{\infty}1\wedge g_2(z)dz \le K_{r,p}( \varsigma^{r/p}+\eta_r). $
At last,
$$ \int_0^{\infty}1\wedge g_3(z)dz =\overline{\sigma}^r\int_0^{\infty}1\wedge\left( \sum_{k=0}^{\infty}\exp\left\{-\frac{z^{2/r}\log\log n_k}{4}\right\}\right)dz=K_r\overline{\sigma}^r. $$
Hence, we conclude that
 \begin{align*} \max_N C_{\upCapc}\left[\max_{n\le N} \left(\frac{ (\sum_{i=1}^n X_i)^+}{a_n}\right)^r\right]
 \le &  2^r\sup_m
 \int_0^{\infty} \Capc\left(\max_{0\le k\le m}I_k\ge z^{1/r}\right)dz \\
 \le &K_{r,p}(\varsigma^{r/p}+\eta_r+\overline{\sigma}^r).
 \end{align*}
 Then, \eqref{eqthMLIL.3} is proved.

 For \eqref{eqthMLIL.4}, notice that \eqref{eq:ad} holds with $p=r$, and so,
 $$ g_2(z)\le c_r z^{-3} \sum_k  \Big(  \frac{n_{k+1}\vSbep[|X|^r]}{a_{n_{k+1}}^r}\Big)^3
 \le    K_r z^{-3}\Big( C_{\upCapc}[|X|^r] \Big)^3.$$
 It follows that
 $ \int_0^{\infty} 1\wedge g_2(z)dz \le K_r C_{\upCapc}[|X|^r] $  by Lemma \ref{lem4}.
 Also,
 $$ \int_0^{\infty} 1\wedge g_1(z)dz\le \sum_{n=1}^{\infty} \int_0^{\infty}\upCapc(|X|\ge z^{1/r}a_n/30) dz
 \le K_r C_{\upCapc}[|X|^r] \sum_{n=1}^{\infty} a_n^{-r}\le K_rC_{\upCapc}[|X|^r].$$
Hence
\begin{align*}
 \max_{N}C_{\upCapc}\left[\max_{1\le n\le N}\left(\frac{ S_n^+ }{a_n}\right)^r \right]
\le   K_r\{C_{\upCapc}[|X|^r]+\overline{\sigma}^r\}\le K_r C_{\upCapc}[|X|^r],
\end{align*}
by noticing that $\overline{\sigma}^r\le  \vSbep[|X|^r] \le C_{\upCapc}[|X|^r]$. The proof of \eqref{eqthMLIL.5} is completed.

 (b) Suppose  \eqref{eqthMLIL.5}. Denote the value on the left hand of \eqref{eqthMLIL.5} by $C_0$. Let $f\in C_{b,Lip}(\mathbb R)$ such that $I\{x>1\}\ge f(x)\ge I\{x>2\}$. Then for $x>0$,
 \begin{align*}
 &  \upCapc\left(2^{r+1}\max_{1\le n\le N}\frac{ |S_n|^r }{a_n^r}>x\right)\ge \upCapc\left(2\max_{1\le n\le N}\frac{ |X_n|^r }{a_n^r}>x\right)\\
 =&1-\lowCapc\left(2\max_{1\le n\le N}\frac{ |X_n|^r }{a_n^r}\le x\right)
 \ge   1-\cSbep\left[\prod_{n=1}^N\left[1-f\left(2\frac{ |X_n|^r }{xa_n^r}\right)\right]\right]\\
 =& 1- \prod_{n=1}^N\cSbep\left[1-f\left(2\frac{ |X_n|^r }{xa_n^r}\right)\right]
 = 1- \prod_{n=1}^N \left[1-\Sbep\left[f\left(2\frac{ |X|^r }{xa_n^r}\right)\right]\right]\\
 \ge & 1-\exp\left\{-\sum_{n=1}^N \upCapc(|X|^r>xa_n^r)\right\}
 \ge   1-\exp\left\{-\int_{16}^{N} \upCapc(|X|^r>xa_y^r)dy\right\}.
 \end{align*}
 Hence
 $$ 2^{r+1}C_0= \max_N\int_0^{\infty} \upCapc\left(2^{r+1}\max_{1\le n\le N}\frac{ |S_n|^r }{a_n^r}>x\right) dx \ge \int_0^{\infty}\big[1-e^{-f(x)}\big]dx, $$
  where   $f(x)=\int_{16}^{\infty} \upCapc(|X|^r>xa_y^r)dy$ is a non-increasing function of $x>0$. Thus, there exists an $x_0>0$ such that $f(x_0)<\infty$. Notice that $(1-e^{-y})/y$ is a non-increasing function of $y>0$. We have that
 \begin{align*}
 & 2^{r+1}C_0\ge   \int_{x_0}^{\infty}[1-e^{-f(x)}]dx\ge \frac{1-e^{-f(x_0)}}{f(x_0)}\int_{x_0}^{\infty}f(x)dx\\
  =& c\int_{1}^{\infty}\int_{16}^{\infty} \upCapc(|X|^r>x_0x a_y^r)dydx=c\int_{16}^{\infty}  \frac{\int_{a_v}^{\infty}ru^{r-1}\upCapc(|X/x_0^{1/r}| > u )du}{a_v^r}dv, \end{align*}
  which, together with \eqref{eqprooflem3.5}, implies   $\eta_{X/x_0^{1/r},r}<\infty$. Hence, $\eta_{X,r}<\infty$.
  Furthermore,
 $$ \max_N\upCapc\left(\max_{1\le n\le N}\frac{|S_n|}{a_n}\ge x\right)\le \frac{C_0}{x^r}<1 \text{ for } x>C_0^{1/r}. $$
By Theorem 5.4 of Zhang (2021),  $\overline{\sigma}_X^2<\infty$,
   $\breve{\mathbb E}[X]=\breve{\mathbb E}[-X]=0$.
   \hfill  $\Box$

%%%%%%%%%%%%%%%%%%%

\end{document}